\documentclass{elsarticle}

\usepackage{lineno,hyperref}
\modulolinenumbers[1]

\usepackage{amssymb}

\journal{Finite Fields and Their Applications}

\newcommand{\Tr}{\textrm{\rm Tr}}

\newcommand{\ord}{\textrm{\rm ord}}
\newcommand{\qord}{\textrm{\rm qord}}

\begin{document}

\def\bbbr{{\rm I\!R}} 
\def\bbbm{{\rm I\!M}}
\def\bbbn{{\rm I\!N}} 
\def\bbbf{{\rm I\!F}}
\def\bbbh{{\rm I\!H}}
\def\bbbk{{\rm I\!K}}
\def\bbbp{{\rm I\!P}}
\def\bbbz{{\mathchoice {\hbox{$\sf\textstyle Z\kern-0.4em Z$}}
{\hbox{$\sf\textstyle Z\kern-0.4em Z$}}
{\hbox{$\sf\scriptstyle Z\kern-0.3em Z$}}
{\hbox{$\sf\scriptscriptstyle Z\kern-0.2em Z$}}}}

\newtheorem{definition}{Definition}
\newtheorem{proposition}{Proposition}
\newtheorem{theorem}{Theorem}
\newtheorem{lemma}{Lemma}
\newtheorem{corollary}{Corollary}
\newtheorem{remark}{Remark}
\newtheorem{example}{Example}
\newtheorem{acknowledgments}{Acknowledgments}

\newenvironment{proof}{\begin{trivlist}\item[]{\em Proof: }}%
{\samepage \hfill{\hbox{\rlap{$\sqcap$}$\sqcup$}}\end{trivlist}}

\begin{frontmatter}

\title{A characterization and an explicit description of all primitive polynomials of degree two over finite fields\tnoteref{mytitlenote}}
\tnotetext[mytitlenote]{Manuscript partially supported by PAPIIT-UNAM IN107423.}


\author{Gerardo Vega}

\address{Direcci\'on General de C\'omputo y de Tecnolog\'{\i}as de Informaci\'on y Comunicaci\'on, Uni\-ver\-si\-dad Nacional Aut\'onoma de M\'exico, 04510 Ciudad de M\'exico, MEXICO \\
(e-mail: gerardov@unam.mx).}

\begin{abstract} 
For polynomials of degree two over finite fields, we present an improvement of Fitzgerald’s characterization (Finite Fields Appl. 9(1):117–121, 2003). We then use this new characterization to obtain an explicit, complete, and simple description of all primitive polynomials of degree two over finite fields.     
\end{abstract}

\begin{keyword}
Primitive polynomials; Order of a polynomial; Cyclotomic classes.
\end{keyword}

\end{frontmatter}


\section{Notation, definitions, and already-known results}\label{secuno}

\noindent
{\bf Notation.} By using $q$ we denote a power of a prime number. From now on, $\gamma$ will denote a primitive element of $\bbbf_{q^2}$ and we fix $\alpha=\gamma^{q+1}$ as primitive element of $\bbbf_{q}$. For any integer $0 \leq i < q-1$, we define ${\cal C}_i^{(q-1,q^2)}:=\gamma^i \langle \gamma^{q-1} \rangle$, where $\langle \gamma^{q-1} \rangle$ denotes the subgroup of $\bbbf_{q^2}^*$ generated by $\gamma^{q-1}$. The $q-1$ cosets, ${\cal C}_i^{(q-1,q^2)}$, are called the {\em cyclotomic classes} of order $q-1$ in $\bbbf_{q^2}$.  

\begin{definition}\label{defcero}
Let $f(x) \in \bbbf_q [x]$ be a polynomial of positive degree with $f(0) \neq 0$. The least positive integer $e$ for which $f(x)$ divides $x^e-1$ is called the {\em order} of $f(x)$ and denoted by \ord$(f(x))$. 
\end{definition}

\begin{definition}\label{defuno}
Let $f(x) \in \bbbf_q [x]$ be a polynomial of positive degree with $f(0) \neq 0$. The least positive integer $\rho$ for which $x^{\rho}$ is congruent modulo $f(x)$, to some element $a \in \bbbf_q^*$, is called the {\em quasi-order} of $f(x)$ and denoted by \qord$(f(x))$. That is $x^{\rho} \equiv a \pmod{f(x)}$ and $\rho=\mbox{\qord}(f(x))$. 
\end{definition}

The following is a quite useful characterization of primitive polynomials over finite fields.

\begin{theorem}\label{TeoNed}
\cite[Theorem 3.18]{Lidl} The monic polynomial $f(x) \in \bbbf_{q}[x]$ of degree $m\geq 1$ is a primitive polynomial over $\bbbf_{q}$ iff $(-1)^m f(0)$ is a primitive element of $\bbbf_{q}$, $\rho=\mbox{\qord}(f(x))=(q^m-1)/(q-1)$, and $x^\rho \equiv (-1)^m f(0) \pmod{f(x)}$.
\end{theorem}

\begin{remark}\label{remuno} 
Let $f(x)$, $\rho$ and $a$ be as in Definition \ref{defuno} and let $e,y,r$ be positive integers such that $e=\rho y + r$, with $0\leq r < \rho$. If $r\neq 0$, then $x^r \equiv g(x) \pmod{f(x)}$, for some non-constant polynomial $g(x) \in \bbbf_q [x]$, and $x^{e}=x^{\rho y}x^r \equiv x^{\rho y}g(x) \equiv a^y g(x) \pmod{f(x)}$. Therefore, for suitable positive integer $e$ and suitable element $d\in \bbbf_q^*$, $x^e \equiv d \pmod{f(x)}$ iff $\rho | e$ and $d=a^{e/\rho}$.
\end{remark}

\section{Preliminary results}\label{secdos}

In the case when the degree of $f(x)$ is two we have:

\begin{proposition}\label{prodos}
Let $f(x) \in \bbbf_q [x]$ be a monic polynomial of degree two and assume for some $a\in \bbbf_q^*$ that $x^{\rho} \equiv a \pmod{f(x)}$, with $\rho=\mbox{\qord}(f(x))$. Let $g(x)=\frac{x^{\rho+1}-ax}{f(x)}=g_{1}x^{\rho-1}+g_{2}x^{\rho-2}+\cdots+g_{\rho-1}x+g_{\rho}$. Then, for $1\leq i\leq \rho-1$, $g_i\neq 0$ and $g_{\rho}=0$  (i.e. apart from the constant term, all terms of $g(x)$ are non-zero). 
\end{proposition}

\begin{proof}
Since $f(x)$ is monic and $g(x)=x\frac{x^{\rho}-a}{f(x)}$, $g_{1}=1$ and $g_{\rho}=0$. Let us assume that $g_k=0$ for some $2\leq k \leq \rho-1$. This means that $f(x) | (x^{k} - g_{k+1})$, which is a contradiction with Definition \ref{defuno} because $k<\rho$.
\end{proof}

\begin{remark}\label{remdos} 
Let $f(x)$, $\rho$ and $g(x)$ be as in Proposition \ref{prodos} and let $e\geq 2$ be an integer. Let $h(x), r(x) \in \bbbf_q [x]$ be the uniquely determined polynomials such that $h(x)=\frac{x^{e}-r(x)}{f(x)}=h_{1}x^{e-2}+h_{2}x^{e-3}+\cdots+h_{e-2}x+h_{e-1}$, with $\deg(r(x))<2$. So note that if $e \geq \rho+1$ then $h_i=g_i$ for $1\leq i\leq \rho$ (in particular $h_{\rho}=g_{\rho}=0$). Alternatively, $e < \rho+1$ iff all the $e-1$ terms of $h(x)$ are non-zero.
\end{remark}

The following result is the key to obtain an explicit and simple description of all primitive polynomials of degree two over finite fields.

\begin{proposition}\label{protres}
Let $f_1(x)=x^2+bx+c, f_2(x)=x^2+m^{j}bx+m^{2j}c \in \bbbf_q [x]$, where $j$ is a positive integer and $b,c,m \in \bbbf_q^*$. Let $e\geq 2$ be an integer and suppose that 

\begin{equation}\label{eqCuno}
x^e \equiv kx+l \pmod{f_1(x)} \;,
\end{equation}

\noindent
for elements $k,l \in \bbbf_q$. Then 

\begin{equation}\label{eqCdos}
x^e \equiv m^{(e-1)j}kx+m^{ej}l \pmod{f_2(x)} \;.
\end{equation}
\end{proposition}

\begin{proof}
Clearly $x^2 \equiv -bx-c \pmod{f_1(x)}$ and $x^2 \equiv m^j(-b)x+m^{2j}(-c) \pmod{f_2(x)}$. Suppose that (\ref{eqCuno}) and (\ref{eqCdos}) hold for $e\geq 2$. Thus, on the one hand, 

$$x^{e+1} \equiv kx^2+lx \equiv k(-bx-c)+lx \equiv (-kb+l)x-kc \pmod{f_1(x)} \;.$$

\noindent
On the other hand, 

\begin{eqnarray}
x^{e+1} &\equiv& m^{(e-1)j}kx^2+m^{ej}lx \nonumber \\
&\equiv& m^{(e-1)j}k(m^j(-b)x+m^{2j}(-c))+m^{ej}lx  \nonumber \\
&\equiv& m^{ej}(-kb+l)x+m^{(e+1)j}(-kc) \pmod{f_2(x)} \;. \nonumber
\end{eqnarray}

\noindent
The result now follows by induction on $e$.
\end{proof}

\begin{corollary}\label{coruno}
Let $f_1(x)$ and $f_2(x)$ be as in Proposition \ref{protres}. Then \qord$(f_1(x))=\mbox{\qord}(f_2(x))$.
\end{corollary}

\begin{proof}
Note that $m^{(e-1)j}\neq 0$. Thus, due to Definition \ref{defuno} and Proposition \ref{protres}, the result follows from the fact that $x^e \equiv l \pmod{f_1(x)}$ iff $x^e \equiv m^{ej}l \pmod{f_2(x)}$, for suitable integer $e\geq2$ and suitable element $l\in \bbbf_q^*$.
\end{proof}

\section{The characterization}\label{sectres}

Through the following result, we identify all irreducible polynomials of degree two that have the same constant term.

\begin{proposition}\label{procinco}
Any irreducible polynomial $x^2+bx+c \in \bbbf_q[x]$ divides $x^{q+1}-c$, and if $q$ is odd, then $x^2-bx+c$ is also a divisor of $x^{q+1}-c$. If $x^2+bx+c$ is primitive then $b\neq 0$ and $c$ is a primitive element of $\bbbf_q$. Finally, if $q$ is odd, then $x^2+bx+c$ is primitive iff $x^2-bx+c$ is primitive.
\end{proposition}

\begin{proof}
Let $\gamma$ and $\alpha$ be as before. Then $c=\alpha^i=\gamma^{(q+1)i}$, for some $0\leq i<q-1$. If $r=\gamma^{i+j(q-1)} \in {\cal C}_i^{(q-1,q^2)}$, for some $0\leq j<q+1$, then $r^{q+1}=\gamma^{i(q+1)}\gamma^{j(q^2-1)}=c$. Since $|{\cal C}_i^{(q-1,q^2)}|=q+1$,

$$x^{q+1}-c=\prod_{r \in {\cal C}_i^{(q-1,q^2)}}(x-r)\;.$$

\noindent
In addition, 

$$r^q=\gamma^{iq+jq^2-jq}=\gamma^{i+i(q-1)-j(q-1)}=\gamma^{i+(i-j)(q-1)} \in {\cal C}_i^{(q-1,q^2)}\;.$$

\noindent
Therefore if $r \in {\cal C}_i^{(q-1,q^2)} \setminus \bbbf_q$, then the irreducible polynomial over $\bbbf_q$, 

$$x^2-\Tr_{\bbbf_{q^2}/\bbbf_q}(r)x+c = (x-r)(x-r^q) \;,$$

\noindent
divides $x^{q+1}-c$, where ``$\Tr_{\bbbf_{q^2}/\bbbf_q}$'' denotes the {\em trace mapping} from $\bbbf_{q^2}$ to $\bbbf_q$. Since $r \in {\cal C}_i^{(q-1,q^2)}$ iff $r^{q+1}=\gamma^{(q+1)i}=c$, there must be an $r' \in {\cal C}_i^{(q-1,q^2)}$ such that $\Tr_{\bbbf_{q^2}/\bbbf_q}(r')=-b$.

Assume that $x^2+bx+c$ is primitive. Note that the quasi-order of any polynomial of the form $x^2-c$, is two. Hence, since $q+1>2$ and due to Theorem \ref{TeoNed}, $b\neq 0$, $c=\alpha^i$ is a primitive element of $\bbbf_q$, and $\gcd(i,q-1)=1$. 

Finally, if $q$ is odd and if $r$ is a root of the primitive polynomial $x^2+bx+c$, then $-r$ is a root of $x^2-bx+c$. Since $r=\gamma^{i+j(q-1)} \in {\cal C}_i^{(q-1,q^2)}$, $-r=\gamma^{\frac{q^2-1}{2}}\gamma^{i+j(q-1)}=\gamma^{i+(j+\frac{q+1}{2})(q-1)}\in {\cal C}_i^{(q-1,q^2)}$. But $\gcd(i,q-1)=1$, so we have $\gcd(i+j(q-1),q^2-1)=1$ iff $\gcd(i+j(q-1)+\frac{q^2-1}{2},q^2-1)=1$, which means that both $r$ and $-r$ are primitive elements of $\bbbf_{q^2}$, which in turn means that both $x^2+bx+c$ and $x^2-bx+c$ are primitive.
\end{proof}

As a first approach, the following result identifies those finite fields for which is quite easy to characterize primitive polynomials of degree two.

\begin{theorem}\label{MiCar1}
Let $f(x)=x^2+bx+c \in \bbbf_{q}[x]$ be irreducible. Suppose that one of the following three conditions holds:

\begin{enumerate}
\item[{\rm (a)}] $q+1=\pi$, for some prime $\pi>2$,
\item[{\rm (b)}] $q+1=2^t$, for some integer $t$, and
\item[{\rm (c)}] $q+1=2 \pi$, for some prime $\pi>2$.
\end{enumerate}

\noindent
Then $f(x)$ is primitive iff $b\neq 0$ and $c$ is a primitive element of $\bbbf_q$.
\end{theorem}

\begin{proof}
If $f(x)$ is primitive, then by Proposition \ref{procinco}, $b\neq 0$ and $c$ is a primitive element of $\bbbf_q$. 

If Condition {\rm (a)} holds, then $q$ is even. Let $r$ be a root of $f(x)$. Then, according to the proof of Proposition \ref{procinco}, $r=\gamma^{i+j(q-1)} \in {\cal C}_i^{(q-1,q^2)}$, for $i$ and $j$ such that $c=\gamma^{(q+1)i}=\alpha^i$. Let $e=i+j(q-1)$ and note that $\gcd(e,q^2-1)>1$ iff $\gcd(e,q^2-1)=\pi=q+1$ iff $\gamma^e \in \bbbf_q^*$. Thus, every element  $r \in {\cal C}_i^{(q-1,q^2)}\setminus \bbbf_q^*$ is a primitive element of $\bbbf_{q^2}$ and $f(x)$ is primitive.

If Condition {\rm (b)} holds, then $q$ is odd. Let $r$ be as before. Since $\gcd(i,q-1)=1$ and $q+1=2^t$, $\gcd(i+j(q-1),q^2-1)=1$. Therefore $r$ is a primitive element of $\bbbf_{q^2}$ and $f(x)$ is primitive. 

Finally, if Condition {\rm (c)} holds, then again $q$ is odd. Let $e=i+j(q-1)$ and note that $\gcd(e,q^2-1)>1$ iff $\gcd(e,q^2-1)=\pi=\frac{q+1}{2}$. But $\pi$ and $i$ are odd, so $4|(q+1+2i)$. If we take $j=\frac{q+1+2i}{4}$, we have

$$e=i+\frac{q+1+2i}{4}(q-1)=\frac{q+1}{2}(\frac{q-1}{2}+i)=\frac{1}{2}[\frac{q^2-1}{2}+(q+1)i]\;.$$

\noindent
Thus, the element $\gamma^e \in {\cal C}_i^{(q-1,q^2)}$ is the root of a non-primitive polynomial. In fact $\gamma^{2e}=\gamma^{\frac{q^2-1}{2}}\gamma^{(q+1)i}=-\alpha^i=-c$, which means that $\gamma^e$ and $\gamma^{qe}$ are the two roots of the irreducible and non-primitive polynomial $x^2+c$. Thus, any element  $r \in {\cal C}_i^{(q-1,q^2)}\setminus \{\gamma^e,\gamma^{qe}\}$ is a primitive element of $\bbbf_{q^2}$ and $f(x)$ is primitive.
\end{proof}

\begin{example}\label{ejuno}
Let $(q,\alpha)=(13,6)$ and note that $\langle 6 \rangle=\bbbf_{13}^*$. Since $1^2-4(2) \equiv 6^1 \pmod{13}$ and $3^2-4(5) \equiv 6^5 \pmod{13}$, the two polynomials $(x^2+x+2)$ and $(x^2+3x+5)$ are both irreducible ($6^1$ and $6^5$ are non-squares in $\bbbf_{13}$). Since $2 \equiv 6^5 \pmod{13}$ and $5 \equiv 6^9 \pmod{13}$, $2$ is a primitive element in $\bbbf_{13}$ and  $5$ is not. Due to Theorem \ref{MiCar1} and since $q+1=14=2(7)$, the first polynomial is primitive while the second one is not.
\end{example}

We are now in position to present our characterization.

\begin{theorem}\label{MiCar2}
Let $f(x)$ be as in Theorem \ref{MiCar1}. In the case that none of the three conditions in Theorem \ref{MiCar1} is satisfied, let $p$ be the smallest odd prime that divides $(q+1)$ and let $h(x), r(x) \in \bbbf_q [x]$ be the uniquely determined polynomials such that $h(x)=\frac{x^{\frac{q+1}{p}+1}-r(x)}{f(x)}$, with $\deg(r(x))<2$. Then $f(x)$ is primitive iff $b\neq 0$, $c$ is a primitive element of $\bbbf_{q}$ and one of the following two conditions holds true:

\begin{enumerate}
\item[{\rm (A)}] $q+1$ is either of the form $q+1=\pi$, $q+1=2^t$, or $q+1=2 \pi$, for some positive integer $t$ and some prime $\pi>2$.
\item[{\rm (B)}] Otherwise, all the $\frac{q+1}{p}$ terms of $h(x)$ are non-zero.
\end{enumerate}
\end{theorem}

\begin{proof}
Condition {\rm (A)} corresponds to the characterization in Theorem \ref{MiCar1}.

Assume Condition {\rm (B)} holds. Let $\rho=\mbox{\qord}(f(x))$ and note that $x^{q+1} \equiv c \pmod{f(x)}$ (recall Proposition \ref{procinco}). Due to Remark \ref{remuno}, $\rho | (q+1)$. If we assume that $f(x)$ is non-primitive, then, by Definition \ref{defuno} and Theorem \ref{TeoNed}, $1<\rho<q+1$. Let $e=\frac{q+1}{p}+1$ and since Condition {\rm (A)} does not holds, $e \geq \rho+1$. Thus, by Proposition \ref{prodos} and Remark \ref{remdos}, $h_{\rho}=0$, a contradiction! Conversely suppose that $f(x)$ is primitive. Thus, by Theorem \ref{TeoNed}, $\rho=q+1$ and, again by Remark \ref{remdos}, all the $e-1=\frac{q+1}{p}<\rho$ terms of $h(x)$ are non-zero.
\end{proof}

\begin{example}\label{ejdos}
Let $(q,\alpha)=(11,2)$. Note that $\langle 2 \rangle=\bbbf_{11}^*$, $q+1=12=2^2(3)$, $p=3$, and therefore $e=\frac{q+1}{p}+1=5$. It is not difficult to see that $x^{12}-2$ is factored in $\bbbf_{11}$ as a product of $6$ irreducible polynomials of degree two:

$$x^{12}-2=(x^2\pm 2x+2)(x^2\pm 4x+2)(x^2\pm 5x+2)\;.$$

\noindent
It is also easy to check that:

\begin{eqnarray}
\frac{x^{5}-(7x)}{(x^2+2x+2)} &=& (x^3+9x^2+2x+0)\;, \nonumber \\
\frac{x^{5}-(10x+8)}{(x^2+4x+2)} &=& (x^3+7x^2+3x+7)\;,\;\mbox{ and}  \nonumber \\
\frac{x^{5}-(6x+1)}{(x^2+5x+2)} &=& (x^3+6x^2+x+5) \;. \nonumber
\end{eqnarray}

\noindent
Thus, owing to Condition {\rm (B)} of Theorem \ref{MiCar2}, we can be sure that the four polynomials, $(x^2\pm 4x+2)$ and $(x^2\pm 5x+2)$, are primitive while the two polynomials, $(x^2\pm 2x+2)$, are non-primitive.
\end{example}

\begin{example}\label{ejtres}
Let $\bbbf_{8}=\bbbf_2(\alpha)$, with $\alpha^3+\alpha+1=0$. Note that $\langle \alpha \rangle=\bbbf_{8}^*$, $q+1=9=3^2$, $p=3$, and therefore $e=4$. It is not difficult to see that $x^{9}+\alpha$ is factored in $\bbbf_{8}$ as a product of $5$ irreducible polynomials:

$$x^{9}+\alpha=(x+\alpha^4)(x^2\!+\alpha x+\alpha)(x^2\!+\alpha^4x+\alpha)(x^2\!+\alpha^5 x+\alpha)(x^2\!+\alpha^6 x+\alpha)\;.$$

\noindent
It is not difficult to check that:

\begin{eqnarray}
\frac{x^{4}+(\alpha^3 x+\alpha^5)}{(x^2\!+\alpha x+\alpha)} &=& (x^2+\alpha x+\alpha^4)\;, \nonumber \\
\frac{x^{4}+(\alpha^5 x)}{(x^2\!+\alpha^4x+\alpha)} &=& (x^2+\alpha^4 x+0)\;, \nonumber \\
\frac{x^{4}+(\alpha x+\alpha)}{(x^2\!+\alpha^5x+\alpha)} &=& (x^2+\alpha^5 x+1)\;,\;\mbox{ and}  \nonumber \\
\frac{x^{4}+(\alpha^4 x+1)}{(x^2\!+\alpha^6 x+\alpha)} &=& (x^2+\alpha^6 x+\alpha^6) \;. \nonumber
\end{eqnarray}

\noindent
Thus, owing to Condition {\rm (B)} of Theorem \ref{MiCar2}, we can be sure that the three polynomials, $(x^2\!+\alpha x+\alpha)$, $(x^2\!+\alpha^5x+\alpha)$ and $(x^2\!+\alpha^6x+\alpha)$, are primitive while the polynomial, $(x^2\!+\alpha^4x+\alpha)$, is non-primitive.
\end{example}

Note that the characterization in \cite[Theorem 1.1]{Fitzgerald} and the characterization in Theorem \ref{MiCar2} are similar in that both characterizations rely on the evaluation of synthetic division of polynomials. However, the synthetic division in Condition (B) of Theorem \ref{MiCar2} involves a dividend polynomial whose degree ($\frac{q+1}{p}+1$) is much smaller than the degree ($q^2-1$) of the dividend polynomial in \cite[Theorem 1.1]{Fitzgerald}. Note also that Condition (A) of Theorem \ref{MiCar2} identifies those finite fields for which it is not necessary to evaluate a synthetic division to determine whether an irreducible polynomial is primitive.

\section{An explicit description of all primitive polynomials of degree two over finite fields}\label{seccuatro}

By considering the polynomials $f(x)$, $h(x)$ and $r(x)$ as in Theorem \ref{MiCar2}, we introduce the Boolean function, $P(q,b,c)$, given by

$$P(q,b,c):=\left\{ \begin{array}{cl}
		\;1 & \mbox{if $b,c\in \bbbf_q^*$, $c$ is a primitive element of $\bbbf_{q}$, and}  \\
		    & \mbox{$f(x)=x^2+bx+x$ is irreducible satisfying one}  \\ 
		    & \mbox{of the two conditions in Theorem \ref{MiCar2},} \\
		\\
		\;0 & \mbox{ otherwise.}
			\end{array}
\right .$$

\noindent
Obviously, a polynomial $f(x)=x^2+bx+c \in \bbbf_q [x]$ is primitive iff $P(q,b,c)=1$. However, despite this obviousness and as we will see below, this function is useful to obtain an explicit and simple description of all primitive polynomials over $\bbbf_q$. To this end, we must first determine whether a polynomial of the form $x^2+\alpha^i x+\alpha \in \bbbf_{q}[x]$ is irreducible or not, where $0\leq i<q-1$ and $\alpha$ is a primitive element in $\bbbf_{q}$. Thus, if $q$ is even and by defining the set:

$$R_0:=\{\:i\:|\:\alpha^{i}=\alpha^{k+1}+\alpha^{q-k-1}, \mbox{ with } 0\leq k<\frac{q}{2}-1\:\}\;,$$

\noindent
it is not difficult to see that $x^2+\alpha^i x+\alpha=(x+\alpha^{k+1})(x+\alpha^{q-k-1})$, for some $1\leq k<\frac{q}{2}-1$, iff $i\in R_0$. Now, if $q$ is odd, then $x^2+\alpha^i x+\alpha$ is irreducible iff $x^2-\alpha^i x+\alpha$ is irreducible. Since $-1=\alpha^{\frac{q-1}{2}}$, we can therefore assume, WLOG, that $0\leq i<\frac{q-1}{2}$. In this context, we denote $i \!\!\!\!\pmod{\frac{q-1}{2}}$ as the integer $0\leq z<\frac{q-1}{2}$ such that $z \equiv i \!\!\!\!\pmod{\frac{q-1}{2}}$. By the definition of the set:

$$R_1:=\{\:i \!\!\!\pmod{\frac{q-1}{2}}\:|\:\alpha^{i}=\alpha^{k+1}+\alpha^{q-k-1}\;,\mbox{ with } 0\leq k<\lfloor \frac{q-1}{4} \rfloor\:\}\;,$$

\noindent
it is not difficult to see that if $i<\frac{q-1}{2}$, then $x^2+\alpha^i x+\alpha=(x+\alpha^{k+1})(x+\alpha^{q-k-1})$, for some $1\leq k<\lfloor \frac{q-1}{4} \rfloor$, iff $i \in R_1$. Now, if $q$ is even we define:

\begin{eqnarray}
B_{0}&:=&\{\:0,1,2,\cdots,q-2\:\} \setminus R_0\;, \nonumber \\
J_{0}&:=&\{\:j\:|\:\gcd(2j+1,q-1)=1, \mbox{ with } 0\leq j<q-1\:\}\;, \nonumber
\end{eqnarray}

\noindent
while if $q$ is odd we define:

\begin{eqnarray}
B_{1}&:=&\{\:0,1,2,\cdots,\frac{q-1}{2}-1\:\} \setminus R_1\;, \nonumber \\
J_{1}&:=&\{\:\frac{j-1}{2}\:|\:\gcd(j,q-1)=1, \mbox{ with } 0\leq j<q-1\:\}\;. \nonumber
\end{eqnarray}

\noindent
Let $\delta=0$ if $q$ is even and $\delta=1$ otherwise. Thereby, as conclusion of our previous discussion, we see that a polynomial of the form $x^2 + \alpha^i x+\alpha$ is irreducible iff either $i \in B_{0}$, if $\delta=0$, or $i \!\!\!\!\pmod{\frac{q-1}{2}} \in B_{1}$, if $\delta=1$. Furthermore note that $\alpha^{2j+1}$ is a primitive element of $\bbbf_q$ iff $j\in J_{\delta}$, $0\leq j<q-1$. With these conclusions in mind, we can now present an explicit and simple description of all primitive polynomials of degree two over a finite field.

\begin{theorem}\label{MiClas1}
With the previous notation, define the subset $I_{\delta} \subseteq B_{\delta}$ given by 

$$I_{\delta} := \{\:i \in B_{\delta} \:|\: P(q,\alpha^{i},\alpha)=1\:\}\;,$$

\noindent
Then, any primitive polynomial of degree two is of the form:

$$x^2 \pm \alpha^{i+j}x+\alpha^{2j+1}\;,$$

\noindent
where $i\in I_{\delta}$ and $j\in J_{\delta}$. 
\end{theorem}

\begin{proof}
It is enough to prove that a polynomial of the form $x^2 + \alpha^{i+j}x+\alpha^{2j+1}$ is primitive. By hypothesis, we know that such a polynomial is primitive and therefore irreducible if $j=0$. On the contrary, for $j>0$ assume that $x^2 + \alpha^{i+j}x+\alpha^{2j+1}$ is reducible. Let $r_1,r_2 \in \bbbf_q^*$ such that $(x-r_1)(x-r_2)=x^2 + \alpha^{i+j}x+\alpha^{2j+1}$. Thus $-r_1-r_2=\alpha^{i+j}$, $r_1r_2=\alpha^{2j+1}$ and $(x-r_1/\alpha^{j})(x-r_2/\alpha^{j})=x^2 + \alpha^{i}x+\alpha$, a contradiction! Now $x^2 + \alpha^{i}x+\alpha$ is primitive, thus, by Theorem \ref{TeoNed} and Corollary \ref{coruno} (take therein $b=\alpha^i$ and $c=m=\alpha$), we have 

$$\mbox{\qord}(x^2 + \alpha^{i+j}x+\alpha^{2j+1})=\mbox{\qord}(x^2 + \alpha^{i}x+\alpha)=q+1\;.$$

\noindent
But $\alpha^{2j+1}$ is a primitive element of $\bbbf_{q}$, thus, again by Theorem \ref{TeoNed}, the polynomial $x^2 + \alpha^{i+j}x+\alpha^{2j+1}$ is primitive. Conversely, if $f(x)=x^2+bx+c \in \bbbf_q[x]$ is primitive, then, $c$ can be expressed as an odd power of $\alpha$; say $c=\alpha^{2j+1}$ for some $j \in J_{\delta}$. Thus, by taking $\alpha^i=b/\alpha^j$, and since $f(x)$ is irreducible and primitive, $i\in I_{0}$, if $\delta=0$, and $i \!\!\!\pmod{\frac{q-1}{2}} \in I_{1}$, if $\delta=1$. 
\end{proof}

\begin{example}\label{ejcuatro}
Let $\bbbf_{8}=\bbbf_2(\alpha)$, with $\alpha^3+\alpha+1=0$. By considering Example \ref{ejtres}, we have $B_0=\{1,4,5,6\}$ and $I_0=\{1,5,6\}$. Clearly $J_0=\{0,1,2,4,5,6\}$. Thus, owing to Theorem \ref{MiClas1}, the eighteen primitive polynomials over $\bbbf_{8}$ are

$$(x^2 + \alpha^{i+j}x+\alpha^{2j+1})\;,$$

\noindent
where $i\in \{1,5,6\}$ and $j \in \{0,1,2,4,5,6\}$. If $\phi$ denotes the Euler totient function, then note that $\phi(q^2-1)/2=18$.
\end{example}

\begin{example}\label{ejcinco}
Let $\bbbf_{9}=\bbbf_3(\alpha)$, with $\alpha^2+\alpha+2=0$. Note that $\langle \alpha \rangle=\bbbf_{9}^*$. It is not difficult to see that $B_1=\{1, 2\}$, and since $q+1=10=2(5)$, $I_1=B_1$ (see Condition {\rm (A)} of Theorem \ref{MiCar2}). Clearly $J_1=\{0,1,2,3\}$. Thus, owing to Theorem \ref{MiClas1}, the sixteen primitive polynomials over $\bbbf_{9}$ are

$$(x^2 \pm \alpha^{i+j}x+\alpha^{2j+1})\;,$$

\noindent
where $i\in \{1,2\}$ and $j \in \{0,1,2,3\}$. Note that $\phi(q^2-1)/2=16$.
\end{example}

\begin{example}\label{ejseis}
Let $(q,\alpha)=(11,2)$ and note that $\langle 2 \rangle=\bbbf_{11}^*$, $2=\alpha^1$, $4=\alpha^2$ and $5=\alpha^4$. Thus, by considering Example \ref{ejdos}, we have $B_1=\{1, 2, 4\}$ and $I_1=\{2,4\}$. Clearly $J_0=\{0,1,3,4\}$. Thus, owing to Theorem \ref{MiClas1}, the sixteen primitive polynomials over $\bbbf_{11}$ are

$$(x^2 \pm \alpha^{i+j}x+\alpha^{2j+1})\;,$$

\noindent
where $i\in \{2, 4\}$ and $j \in \{0,1,3,4\}$. Again, note that $\phi(q^2-1)/2=16$.
\end{example}

The most difficult part of Theorem \ref{MiClas1}, is the determination of the subset $I_{\delta} \subseteq B_{\delta}$. In the worst case, this subset requires $|B_{\delta}|$ polynomial divisions of the type in Condition {\rm (B)} of Theorem \ref{MiCar2}. Finally, note that the set $B_{\delta}$ is relatively easily determined by the set $R_{\delta}$ using only the arithmetic in $\bbbf_q$.

\end{document}